\documentclass[14pt]{article}

\usepackage{mathtext}
\usepackage[T1,T2A]{fontenc}
\usepackage[utf8]{inputenc}
\usepackage[english]{babel}

\usepackage{amsfonts, amsmath, amssymb, amsthm, amscd}
\usepackage{mathtools}
\usepackage{xcolor}
\usepackage{ulem}
\usepackage{array}

\usepackage{geometry}
\usepackage{graphicx}
\usepackage{hyperref}
\usepackage{cancel}
\hypersetup{pdfstartview=FitH,  linkcolor=blue,urlcolor=urlcolor, citecolor=blue, colorlinks=true}

\theoremstyle{theorem}
\newtheorem{theorem}{Theorem}[section]
\theoremstyle{theorem}

\theoremstyle{definition}

\theoremstyle{lemma}

\newtheorem{corollary}{Corollary}[section]
\newcommand{\Label}[1]{\label{#1}}

\usepackage{fancyhdr}
\newcommand\RR{\mathbb{R}}
\newcommand\CC{\mathbb{C}}
\newcommand{\cH}{\mathcal{H}}
\newcommand{\cF}{\mathcal{F}}
\newcommand{\cJ}{\mathcal{J}}
\newcommand{\cD}{\mathcal{D}}
\newcommand{\pr}{\prime}
\newcommand{\Qp}{\mathbb Q_p}
\newcommand{\Fx}{\cF_{x\to \xi}}
\newcommand{\Fi}{\cF^{-1}_{\xi\to x}}

\newcommand{\wh}[1]{\widehat{#1}}
\newcommand{\col}{\colon}

\newcommand\bn{\,|\!\!|\!\!|\,}

\newcommand{\cK}{\mathcal K}

\newcommand{\cU}{\mathcal U}

\newcommand{\vph}{\varphi}

\newcommand{\ga}{\gamma}
\newcommand{\be}{\beta}

\newcommand{\OO}{\Omega}

\newcommand{\al}{\alpha}

\newcommand{\bs}{\backslash}

\numberwithin{equation}{section}

\begin{document}

\title{Non-Archimedean Kelvin Transformation}

\author{
        Alexandra V. Antoniouk$^1$
        \and
        Anatoly N. Kochubei$^2$
 }

\author{Alexandra V. Antoniouk\\
Institute of Mathematics, National Academy of Sciences of Ukraine,\\ Tereshchenkivska 3, Kyiv, 01024 Ukraine\\
American University Kyiv, Poshtova Sq 3, 04070, Kyiv, Ukraine\\
E-mail: antoniouk.a@gmail.com
\and
Anatoly N. Kochubei\\
Institute of Mathematics, National Academy of Sciences of Ukraine,\\ Tereshchenkivska 3, Kyiv, 01024 Ukraine \\
E-mail: kochubei@imath.kiev.ua}

\date{}

\maketitle

\begin{abstract}
We introduce and study an analog of the Kelvin transformation connected with the Vladimirov-Taibleson operator acting on real- or complex-valued functions on a space $K^n$ over a non-Archimedean local field $K$.
\end{abstract} 


\section{Introduction} \label{sec:1}
The classical Kelvin transform is an important tool of mathematical analysis. In the framework of $\al$-fractional potential theory \cite{La,RS}, the Kelvin transform if a function $u$ on $\RR^n$, $n\geq 2$ is
\begin{equation}\Label{1-1}\
\big(\cK u\big)(y)= |y|^{\al - n}u(Ty),
\end{equation}
where $0< \al\leq 2$, $T$ is the inversion with respect to the unit sphere defined as
\begin{equation}\Label{1-2}\
Tx=\frac{x}{|x|^2}.
\end{equation}

The case where $\al = 2$ is related to the theory of harmonic functions (see, for example, \cite{ABR}). If $\OO$ is an open subset of $\RR^n\bs \{0\}$ and $\OO^*=\{Tx\col x\in \OO\}$, then a function $u$ is harmonic on $\OO$ if and only if $\cK u$ is harmonic on $\OO^*$. This property has a number of applications. It is extended to $\al$-harmonic functions related to the fractional Laplacians and other operators. There is a deep theory devoted to probabilistic interpretation of the Kelvin transform \cite{BZ, ACGZ}.

In this note we introduce a kind of Kelvin transformation connected to the Vladimirov-Taibleson operator acting on real- or complex-valued functions on a space $K^n$ over a non-Archimedean local field $K$:

\begin{equation}\Label{1-3}\
\big(D^{\al.n}u\big) (x)= c_{n,\al}\int\limits_{K^n} \frac{u(x)-u(y)}{\|x-y\|_{K^n}^{n+\al}}\, dy, \quad x\in\OO\subset K^n,
\end{equation}
where $n\geq 2$, $\Vert z\Vert_{K^n}=\max\limits_{1\leq j \leq n}|z_j|_K$  for $z=(z_1,\ldots,z_n)\in K^n$, $|\cdot|_K$ is the normalized absolute value, $dy$ is the Haar measure,
\[c_{n,\al} = \frac{q-1}{1-q^{-\al-n}}, \quad \al > 0,
\]
$q$ is the residue cardinality of $K$, and the precise meaning of the inegral in \eqref{1-3} depends on smoothness assumptions regarding the function $u$. A model example of the field $K$ is the field $\Qp$ of $p$-adic numbers.

The operator \eqref{1-3} and equations containing it have been investigated by many authors. Its theory is to some extent parallel to that of fractional Laplacians of real analysis, though the details are often very different; see the books \cite{VVZ, K2001, T, AKS, KKZ, Z2016, Z24} and many recent papers, such as \cite{BGPW, AKN, AKSe, K2023} and others.

The formulas \eqref{1-1} and \eqref{1-2} cannot be reproduced in the non-Archimedean situation, since an element $x\in K^n$ cannot be divided in a straightforward way by a real number $|x|_K^2$. An interpretation of the latter as a $p$-adic number does not help: $p$-adic geometry is quite different from the real one.

Fortunately, the non-Archimedean space $K^n$ possesses the property absent in $\RR^n$. On the set $K^n$, there exist various field structures, in particular the structure of an unramified field extension (see Section \ref{sec:2}). This makes it possible to define an analogue of the Kelvin transformation, connected with the notion of $\al$-harmonic functions, solutions fo the equation $\big(D^{\al,n}u\big)(x)=0$ on subsets of $K^n$.

\section{Local fields} \label{sec:2}

A (non-Archimedean) local field is non-discrete totally disconnected locally compact topological field, Such a field $K$ is isomorphic to a finite algebraic extension of the field $\Qp$ of $p$-adic numbers, if $K$ has characteristic $0$, or ot the field of formal Laurent series with coefficients from a finite field, if $K$ has a positive characteristic. See \cite{T, S1979, W1967, K2001, FV} for the basic notions and results regarding the local fields.

On any local field there exists an absolute value $|\cdot|_K$, such that: 1) $|x|_K=0$ if and only if $x=0$; 2) $|xy|_K=|x|_K|y|_K$; 3) $|x+y|_K\leq \max\{|x|_K, |y|_K\}$. This last inequality, called ultrametric, implies that $|x+y|_K=|x|_K$, if $|y|_K<|x|_K$.

The ring $O=\{x\in K\col |x|_K\leq 1\}$ is called the ring of integers of $K$. The ideal $P=\{x\in K\col |x|_K <1\}$ contains such an element $\be$ that $P=\be O$. The quotient ring $\bar K=O/P$ is actually a finite field called the residue field. The absolute value is called normalized, if $|\be|_K=q^{-1}$ where $q$ is the cardinality of $O/P$.

If a local field $K$ is a subfield of a local field $L$, then $L$ is called an extension of $K$, which is denoted $L/K$. The extension $L$ can be considered as a vector space over $K$. If $L$ is finite-dimensional over $K$, then $L/K$ is called the finite extension of degree $n$, where $n=\text{dim}\  L$ over $K$. Any basis of $L$ over $K$ is called a basis of the extension.

An operator of multiplication in $L$ by an element $\xi$ can be considered as a linear operator on a $K$-vector space. If the linear function $\xi \mapsto Tr(\xi)$ does not vanish identically, then the extension is called separable. All finite extensions of the field of characteristic zero are separable. On the other hand, the notion of separability makes sense also for the finite fields $\bar K, \bar L$.

There exists a detailed theory of extension of local fields; see \cite{FV, W1967}. Here we will need only the unramified extensions.

A finite extension $L$ of a local field $K$ is called unramified, if $\bar L/\bar K$ is a separable extension of the same degree, as $L/K$. Any local field $K$ has a unique (up to the isomorphism) unramified extension of each given degree $n\geq 2$. The cardinality of the residue field $\bar L$ in this extension equals $q^n$ where $q$ is the residue cardinality $\text{card}\, (O/P)$ in $K$. Each unramified extension is separable.

The unramified extension of degree $n$ has, as a $K$-vector space, a canonical basis consisting of representatives of a basis in $\bar L$ over $\bar K$. Let $x\in L$ have the coefficients $x_1,\ldots,x_n \in K$ of its expansion with respect to the canonical basis. Then the normalized absolute value $|x|_L$ has the representation \cite{T1976}:
\begin{equation}\Label{28}\
|x|_L=\Big(\max\limits_{1\leq j\leq n} |x_j|_K\Big)^n.
\end{equation}
Another absolute value on $L$, $\| x\|_L=|x|^{1/n}_K\equiv\max\limits_{1\leq j \leq n}|x_j|_K$, extends the absolute value from $K$. Comparing the expression for $\|x\|_L$ and $\|x\|_{K^n}$, we see that the above canonical basis is an orthogonal basis (in the non-Archimedean sense \cite{Sc1984} in $L$ considered as Banach space with the norm $\|\cdot \|_L$. Therefore (see also \cite{K2021}) the canonical basis defines an isometric linear isomorphism between $L$ and $K^n$.

In particular, this isomorphism defines an invariant measure on the additive group $K^n$, and the measure of the unit ball in $K^n$ equals 1. By the uniqueness of the Haar measure, $dx = dx_1\ldots dx_n$. Moreover, this isomorphism identifies the operator $D^{\al,n}$ defined in \eqref{1-3} with the operator
\begin{equation}\Label{2-1}\
\Big(D_L^\ga u\Big)(x)=\frac{1-q^{n\ga}}{1-q^{-n(\ga+1)}}\int\limits_L \frac{u(z)-u(x)}{|z-x|_L^{\ga+1}}\, dz, \quad \ga=\frac{\al}{n}.
\end{equation}
Thus\footnote{Why in formule \eqref{2-1} it appears $n$ in the constant in front of the inegral? I thought that the operator $D_L^\be$ is like one-dimensional opeator, with the corresponding normanizing constant?}, a multi-dimensional operator over $K$ is equivalent to a one-dimensional operator over $K$ is equivalent to a one-dimensional operator over a bigger field. Of course, this is a purely non-Archimedean phenomenon.

\section{Fourier analysis and Sobolev-type spaces} \label{sec:3}

The additive group of a local field $K$ is self-dual, so that the Fourier analysis on $K$ is similar to the classical one. The details may be found in particular in \cite{K2001, KKZ, VVZ, Z24}.

Let $\chi$ be a fixed non-constant additive character. We assume that $\chi$ had rank zero, so that $\chi(x)\equiv 1$ on $O$, and $\chi(x_0)\neq 1$ for some $x_0\in K$ with $|x_0|_K=q$.

The Fourier transform of a complex-valued function $f\in L^1(K)$ is defined as

\[\big(\Fx f\big)(\xi)= \widehat{f}(\xi)=\int\limits_{K}\chi(x\xi)f(x)\, dx, \quad \xi \in K.\]

If $\hat{f}\in L^1(K)$, then we have an inversion rule
\[f(x)=\int\limits_K \chi (-x\xi)\widehat{f}(\xi)\, d\xi.\]

The Fourier transform preserves the Bruhat-Schwartz space $\cD(K)$ of test functions, consisting of locally constant functions with compact supports. The local constancy of a function $f\col K\mapsto \CC$ means the existence of such an integer $k$ that for any $x\in K$
\[f(x+x^\pr)=f(x), \quad \text{whenever}\quad |x^\pr|_K\leq q^{-k}.\]

The above notion extend naturally to functions $K^n\mapsto \CC$.

For $\vph\in\cD(K^{{n}})$,
\[\big(D^{\al,n}\vph\big)(x)=\Fi\big(\|\xi\|^\al_{K^n}\Fx\vph\big)(x).\] We will use a consruction of Sobolev-type spaces of complex-valued functions suggested by Zunigo-Galindo (for the details see \cite{KKZ, Z24}.

Let $C_0(K^n)$ be the space of all continuous functions on $K^n$ tending to zero at infinity by the filter of complements to compact subsets. For $\vph, \psi\in \cD(K^n)$, define the scalar product
\begin{equation}\Label{3-1}\
\langle \vph, \psi \rangle_\ell=\int\limits_{K^n}\bn\xi\bn_{K^n}^\ell \wh{\vph}(\xi)\overline{\wh{\psi}(\xi)}\, d\xi, \quad \ell = 0,1,2, \ldots,
\end{equation}
where $\bn\xi\bn_{K^n} =\max (1,\|\xi\|_{K^n})$. Let $\cH_\ell=\cH_\ell(K^n)$ be the completion of $\cD(K^{{n}})$ with respect to the scalar product \eqref{3-1}. We have continuous embedding $\cH_m \hookrightarrow \cH_\ell$ for $\ell \leq m$. Next, we set
\[\cH_\infty = \cH_\infty (K^n)=\bigcap_{\ell = 0}^\infty\cH_\ell.\]

We have $\cH_0 = L^2$, $\cH_\infty \subset L^2$.

The scalar products \eqref{3-1} define a system of seminorms on $\cD(K^n)$, so that $\cD(K)$ becomes a metric space. Its completion coincides with  $\cH_\infty (K^n)$. Thus $\cH_\infty (K^n)$ can be considered as a locally convex vector space continuously imbedded in $C_0(K^n)$. The mapping $D^{\al,n}$ is a bicontinuous  automorphism of $\cH_\infty (K^n)$.

\section{Kelvin transformation}\label{sec:4}

Denote by $\cU\colon K^n \to L$ the above isometric isomorphism. An inversion in $K^n$ is defined as
\begin{equation}\Label{4-1}\
\cJ= \cU^{-1} \Delta\, \cU,
\end{equation}
where $\Delta (x)=\dfrac{1}{x}$, $x\in L\bs \{0\}$. Note that the classical inversion \eqref{1-2} is a conformal mapping, see \cite{ABR}. In $p$-adic mathematical physics the conformal group is defined as the one generated by $\Delta$ and isometries \cite{Le}; in this sense our inversion \eqref{4-1} is a conformal mapping.

We call the transformation
\[\big(\cK u\big)(x)=\| x\|^{\al - n}_{K^n}u(\cJ(x)), \quad 0 \neq x\in K^n\]
the non-Archimedean Kelvin transformation.

Let us consider the operator $D^\ga = D^\ga_L$ ($\ga\neq 1$) of the form \eqref{2-1} corresponding to the field $L$. Let us denote $x^*=\dfrac{1}{x}$ for $x\in L$ and $f^*(x)=f(x^*)$. We calculate the action of operator $\big(D^\ga f^*\,\big)(x^*)$, on the function $f^*$, which is sufficiently regular, so that the calculations make sense. For example, we may assume $f\in\cH_\infty (L)$. The nonlinear change of variables in an integral is substantiated in \cite[Appendix~A7]{Sc1984}.

First we consider the operator
\[\big(D^{-\ga} f^*\,\big)(x^*)=d_{n,\ga}\int\limits_L|x^*-y|_L^{\ga -1} f\big(\frac{1}{y}\big)\, dy, \quad x\neq 0,\]
where
 \[d_{n,\ga}=\frac{1-q^{-n\ga}}{1-q^{n(\ga-1)}}.\]
 The operator $D^{-\ga}$ is an inverse operator to $D^{\ga}$ \cite{VVZ}. We have

\begin{align}\Label{29-2a}\ \nonumber
&\big(D^{-\ga} f^*\big)(x^*)=d_{n,\be}\int\limits_L\Big|x^*-\frac{1}{z}\Big|_L^{\ga-1}\Big|\frac{1}{z^2}\Big|_Lf(z)\,dz=\\
\nonumber
&= d_{n,\ga}\int\limits_L |z|_L^{-\ga-1}|x^* z-1|_L^{\ga -1}f(z)\, dz=\\
\nonumber
&=d_{n,\ga}|x^*|_L^{\ga-1}\int\limits_L|z|_L^{-\ga-1}|z-x|_L^{\ga-1} f(z)\,dz =\\
&= |x^*|_L^{\ga-1} D^{-\ga}\Big(|\cdot|_L^{-\ga-1}f(\cdot)\Big)(x).
\end{align}

If $f=D^\ga u$, then
 \[u\,\Big(\frac{1}{x}\,\Big)=|x|_L^{1-\ga} D^{-\ga}\Big(|\cdot|_L^{-\ga-1}D^\ga u(\cdot)\Big)(x),
 \]
so that
\[ D^{\ga}\Big( |x|_L^{\ga-1}u(x^*)\Big)(x)= |x|_L^{-\ga-1}D^\ga_x u(x),
 \]
or, equivalently,
\begin{equation}\Label{4-1a}\
\big(D^{\ga} u\big)(x)= |x|_L^{\ga+1}D^\ga\Big( |x|_L^{\ga-1} u(x^*)\Big)(x).
 \end{equation}

 The identity \eqref{4-1a} can be interpreted for functions on $K^n$, if we set $\ga=\dfrac{\al}{n}$ and identify $L$ with $K^n$ using the expansion with respect to a canonical basis. Thus, we have proved the following result.

 \begin{theorem}
 For any $u\in \cH_\infty (K^n),$
 \[\Big( D^{\al,n}u\Big) (x)=\|x\|^{\al + n}D^{\al,n}\big(\cK u\big)(x)\]
 \end{theorem}

 As usual we call a function $u\in\cH_\infty(K^n)$ $\al$-harminuc on a subset $G\subset K^n$, if $\big( D^{\al,n} u\big)(x)=0$ for all $x\in G$. Note that an $\al$-harmonic function (on a subset) must be defined on the whole space $K^n$. For explicit example of such function see the eigenfunctions of $D^{\al,n}$ \cite{VVZ, K2001, KKZ, BGPW}.

 \begin{corollary} Let $u$ be $\al$-harmonic function on $G$, $G\ni 0$. Then the function $\cK$ is $\al$-harmonic on the set $\{\cJ (x)\col x\in G\}$.
 \end{corollary}

\bigskip
\section*{Acknowledgements}
The first author acknowledges the funding support in the framework of the project ``Spectral Optimization: From Mathematics to Physics and Advanced Technology'' (SOMPATY) received from the European Union’s Horizon 2020 research and innovation programme under the Marie Skłodowska-Curie grant agreement No 873071. The second author acknowledges the financial support by the National Research Foundation of Ukraine (Project number 2023.03/0002).


\end{document}